\documentclass[showpacs,preprintnumbers,amsmath,amssymb]{revtex4}

\usepackage{graphicx}
\usepackage{dcolumn}
\usepackage{bm}
\usepackage{mathrsfs}
\usepackage{amsmath}
\usepackage{enumerate}
\usepackage{color}
\newtheorem{theorem}{Theorem}[section]

\newenvironment{proof}[1][Proof.]{\begin{trivlist}
     \item[\hskip \labelsep {\bfseries #1}]}{\end{trivlist}}

\DeclareMathOperator*{\foo}{K}

\begin{document}


\title{The mathematical foundations of the asymptotic iteration method}

\author{D. Batic}
\email{davide.batic@ku.ac.ae}
\affiliation{%
Department of Mathematics,\\  Khalifa University of Science and Technology,\\ Main Campus, Abu Dhabi,\\ United Arab Emirates}
\author{M. Nowakowski}
\email{mnowakos@uniandes.edu.co}
\affiliation{
Departamento de Fisica,\\ Universidad de los Andes, Cra.1E
No.18A-10, Bogota, Colombia
}

\date{\today}

\begin{abstract}
We introduce a new approach to the the asymptotic iteration method (AIM) by means of which we establish the standard AIM connection with the continued fractions technique and we develop a novel termination condition in terms of the approximants. With the help of this alternative  termination condition and certain properties of continuous fractions, we derive a closed formula for the asymptotic function $\alpha$ of the AIM technique in terms of an infinite series. Furthermore, we show that such a series converges pointwise to $\alpha$ which, in turn, can be interpreted as a specific term of the minimal solution of a certain recurrence relation. We also investigate some conditions ensuring the existence of a minimal solution and hence, of the function $\alpha$ itself.
\end{abstract}

\maketitle

\section{Introduction}\label{Intr}
\cite{Ciftci} developed the so-called asymptotic iteration method (AIM) in order to solve certain eigenvalues problems described by second order linear ODEs of the form
\begin{equation}\label{ode1}
y^{''}=\lambda_0(x)y^{'}+s_0(x)y,\quad \lambda_0,s_0\in C^\infty (I),\quad I\subseteq\mathbb{R}
\end{equation}
with $\lambda_0\neq 0$ for any $x\in I$ and furthermore, $\lambda_0$ and $s_0$ do not need to be bounded at the endpoints of the interval $I$. The spectral parameter is usually encoded in the function $s_0$. The smoothness condition on $\lambda_0$, and $s_0$ is not too restrictive because in most applications such functions are usually represented by some rational functions which can be always made smooth by suitably choosing their domain of definition. Moreover, \cite{Ciftci} tacitly assumed that the general solution of (\ref{ode1}) is also smooth on the interval $I$ where both $\lambda_0$ and $s_0$ are smooth. This assumption is reasonable because the solution of (\ref{ode1}) is discontinuous at most at those points where $\lambda_0$ and $s_0$ are not continuous. \cite{Ciftci} differentiates equation (\ref{ode1}) $n$-times with respect to the independent variable to end up with  
\begin{equation}\label{ode2}
y^{(n+2)}=\lambda_n(x)y^{'}+s_n(x)y
\end{equation}
where the functions $\lambda_n$ and $s_n$ are determined by the following coupled system of recurrence relations
\begin{eqnarray}\label{prasec}
\lambda_n&=&\lambda_{n-1}^{'}+\lambda_0\lambda_{n-1}+s_{n-1},\quad\forall n=1,2,\cdots\label{rr1}\\
s_n&=&s^{'}_{n-1}+s_0\lambda_{n-1}.\label{rr2}
\end{eqnarray}
At this point, \cite{Ciftci} considers the ratio of the $(n+2)$-th and $(n+1)$-th derivatives which can be written by means of (\ref{ode2}) as
\begin{equation}\label{2a}
\frac{y^{(n+2)}}{y^{(n+1)}}=\frac{\lambda_n}{\lambda_{n-1}}\frac{y^{'}+\frac{s_n}{\lambda_n}y}{y^{'}+\frac{s_{n-1}}{\lambda_{n-1}}y}
\end{equation}
together with the asymptotic condition in \cite{Ciftci}
\begin{equation}\label{ascond}
\frac{s_n}{\lambda_n}=\frac{s_{n-1}}{\lambda_{n-1}}:=\alpha(x),\quad n\gg 1,
\end{equation}
which is equivalent to {\bf{requiring}} that there exists an $N\in\mathbb{N}$ such that for all $n\geq N$
\begin{equation}\label{condo}
\frac{s_n(x)}{\lambda_n(x)}=\alpha(x)
\end{equation}
for all $x\in I$. Finally, if the condition (\ref{condo}) holds and we make use of (\ref{rr1}) and (\ref{rr2}), we conclude that (\ref{2a}) admits the solution
\begin{equation}\label{y1n}
y^{(n+1)}(x)=C_1\lambda_{n-1}(x)\mbox{exp}\left(\int_x\left[\lambda_0(t)+\alpha(t)\right]~dt\right)
\end{equation}
for all $n\geq N$. Substituting (\ref{y1n}) into $y^{(n+1)}=\lambda_{n-1}y^{'}+s_{n-1}y$ and applying the asymptotic condition (\ref{condo}), we come up with the following representation for the solution {\bf{of}} (\ref{ode1}), namely
\begin{equation}\label{gs}
y(x)=\mbox{exp}\left(-\int_x\alpha(\rho) d\rho\right)\left[C_1\int_x\mbox{exp}\left(\int_\tau\left[\lambda_0(t)+\alpha(t)\right]~dt\right)d\tau+C_2\right].
\end{equation}
\cite{Ciftci} claims that (\ref{gs}) is the general solution of (\ref{ode1}). Such a claim should be taken with some caution. First of all, (\ref{gs}) would need that (\ref{condo}) holds for all $n\in\mathbb{N}$ which is not the case in general. Moreover, if we replace (\ref{gs}) into (\ref{ode1}), we find that (\ref{gs}) is a solution of the original ODE if and only if the unknown function $\alpha$ is a solution of the Riccati equation
\begin{equation}\label{Riccati}
\alpha^{'}=\alpha^2(x)+\lambda_0(x)\alpha(x)-s_0(x).
\end{equation}
It is straightforward to verify that a class of solutions to (\ref{ode1}) can be represented as in (\ref{gs}) with $\alpha$ solution to the Riccati equation (\ref{Riccati}) if and only if the second order differential operator in (\ref{ode1}) admits the factorization
\begin{equation}
\left(\frac{d}{dx}+\alpha(x)\right)\left(\frac{d}{dx}+\beta(x)\right)y=0,\quad
\beta=-\lambda_0(x)-\alpha(x).
\end{equation}
This property is completely independent on the AIM asymptotic condition (\ref{ascond}) whose origin can be related to the requirement that the general solution of (\ref{ode1}) is real-analytic on some interval $I$. As we will see here below, this observation is a key ingredient in order to derive new asymptotic conditions for the AIM method. To see that, let us take any $x_0\in I$. Note that the interval $I$ has been already so defined that $\lambda_0$ never vanishes there.  Consider the initial condition $y(x_0)=\mathfrak{c}_1$. Then, we can set up a Taylor series about $x_0$ for the solution of (\ref{ode1})
\begin{equation}\label{Taylor}
y(x)=\mathfrak{c}_1+(x-x_0)\sum_{n=0}^\infty f_n(x),\quad f_n(x)=\frac{y^{(n+1)}(x_0)}{(n+1)!}(x-x_0)^n
\end{equation}
and
\begin{equation}\label{ratio}
\left|\frac{f_{n+1}(x)}{f_n(x)}\right|=\frac{|x-x_0|}{n+2}\left|\frac{y^{(n+2)}(x_0)}{y^{(n+1)}(x_0)}\right|=\frac{|x-x_0|}{n+2}\left|\frac{\lambda_n(x_0)}{\lambda_{n-1}(x_0)}\right|\left|\frac{y^{'}(x_0)+\mathfrak{c}_1\frac{s_n(x_0)}{\lambda_n(x_0)}}{y^{'}(x_0)+\mathfrak{c}_1\frac{s_{n-1}(x_0)}{\lambda_{n-1}(x_0)}}\right|,
\end{equation}
where in the last step we used (\ref{2a}) evaluated at the point $x_0$ together with the initial condition $y(x_0)=\mathfrak{c}_1$. If the asymptotic condition
\begin{equation}\label{asympt2}
\lim_{n\to\infty}\frac{s_n(x_0)}{\lambda_n(x_0)}=\alpha(x_0),\quad \alpha(x_0)<\infty
\end{equation}
and the asymptotic expansion
\begin{equation}\label{asy}
\left|\frac{\lambda_n(x_0)}{\lambda_{n-1}(x_0)}\right|=\mathfrak{a}_0+\mathfrak{a}_1 n+\mathcal{O}\left(\frac{1}{n^\gamma}\right),\quad\mathfrak{a}_i\geq 0~\forall i=1,2\quad, \gamma>0
\end{equation}
hold, then the Taylor series (\ref{Taylor}) will converge in some open interval about $x_0$ with convergence radius $\rho=1/\mathfrak{a}_1$. At this point a remark is in order. Pointwise convergence of the sequence $(s_n/\lambda_{n})_{n\in\mathbb{N}}$ together with the asymptotic expansion (\ref{asy}) is sufficient to ensure convergence of (\ref{Taylor}) for any $x_0\in I$, and therefore, the smoothness of the solution to (\ref{ode1}) on the interval $I$. However, if we also require that the general solution to (\ref{ode1}) admits the representation (\ref{gs}), we need to demand that $\alpha$ is at least continuous on $I$ but this is equivalent to {\bf{imposing}} uniform convergence for the sequence $(s_n/\lambda_{n})_{n\in\mathbb{N}}$ on the interval $I$. {\bf{Note}} that \cite{Ciftci,Ismail} do not discuss which kind of convergence should be imposed on (\ref{ascond}). Let us introduce the sequence $(\alpha_n)_{n\in\mathbb{N}}$ with $\alpha_n=s_n/\lambda_n$. If such a sequence converges at least pointwise on $I$, then the {\bf{new}} AIM asymptotic condition can be formulated as
\begin{equation}\label{equiv}
\lim_{n\to\infty}\left[\alpha_n(x)-\alpha_{n-1}(x)\right]=0
\end{equation}
or, in the case a termination behaviour occurs as 
\begin{equation}\label{term}
\alpha_N(x)-\alpha_{N-1}(x)=0
\end{equation}
for some positive integer $N$. If the termination condition (\ref{term}) is not satisfied by the particular problem at hand, it is necessary to consider (\ref{equiv}). \cite{Ciftci,Cho1,Cho2} noticed that the stability and number of iterations needed in (\ref{equiv}) to compute eigenvalues of certain eigenvalues problems rely in a sensitive way on the choice of the point $x_0\in I$. This {\bf{signals}} that uniform convergence for the sequence $(\alpha_n)_{n\in\mathbb{N}}$ will be more the exception than the rule. In that regard, we would like to mention that \cite{Cama} showed that in the case (\ref{ode1}) is a generic Fuchsian differential equation in the complex plane, the aforementioned sequence converges uniformly to the logarithmic derivative of the holomorphic solution to (\ref{ode1}) in compact subsets of $U_i=\{z\in\mathbb{C}~|~|z-z_i|<|z-z_j|~|~j=1,2,\cdots,r,~j\neq i\}$ where the index $r$ counts the singular points of (\ref{ode1}). Concerning the optimal choice for $x_0$, \cite{Ismail} gives indeed a condition represented by equation (56) therein, which is mutatis mutandis equivalent to (\ref{equiv}) and specifies when the AIM should work, however no algorithm or method is offered by \cite{Ismail} leading to the optimal choice of $x_0$. Moreover, the connection between AIM and continuous fractions was also studied in \cite{Matamala}. Last but not least, the power of the AIM becomes evident if we recall that the new branch of quasinormal modes for a massless scalar field in the Schwarzschild metric derived analytically in \cite{34,35,36} can also be obtained numerically by using the improved AIM \cite{Ciftci,Cho1,Cho2,37}. The paper is organised as follows: in Section~\ref{Abt2} we prove the uniqueness of the solution to (\ref{prasec}) subject to an initialisation condition. In Section~\ref{Abt3}  we show how the unknown function $\alpha$ appearing in (\ref{Riccati}) can be represented in terms of a certain continued fraction and we also discuss the convergence problem of the latter. We also construct a new formula for the unknown function $\alpha$ represented by (\ref{zve}) and we show that if the conditions in Theorem~\ref{TH3} are met, the sequence with terms $\alpha_n=s_n/\lambda_n$ in the AIM will converge pointwise on some interval $I$ to the function $\alpha$ which in turn can be interpreted as the term $\psi_{-1}(x)$ of the minimal solution $\psi_n(x)$ of the recurrence relation (\ref{3rr}). We also investigate several conditions ensuring the existence of a minimal solution and we observe that they have the effect to shrink the interval over which the optimal value for $x_0$ must be chosen when we use the AIM.

\section{The coupled system of recurrence relations (\ref{rr1}) and (\ref{rr2})}\label{Abt2}
Let $I_2$ denote the $2\times 2$ identity matrix. Define
\begin{equation}
\Phi_n:=\left(
\begin{array}{c}
\lambda_n\\
s_n
\end{array}
\right)~\forall n\in\mathbb{N},\quad T:=I_2\frac{d}{dx}+\mathfrak{A}(x),\quad \mathfrak{A}(x):=\left(
\begin{array}{cc}
\lambda_0 & 1\\
s_0 & 0
\end{array}
\right).
\end{equation}
Note that at this step the differential operator $T$ is a formal linear operator defined on $\mathbb{C}(x)\times\mathbb{C}(x)$ with $\mathbb{C}(x)$ denoting the vector space of all rational functions in the real variable $x$ where the polynomials in the numerator and denominator may have complex coefficients.  Then, the recurrence relations (\ref{rr1}) and (\ref{rr2}) can be written in matrix form as 
\begin{equation}\label{msystem}
\Phi_n=T\Phi_{n-1}\quad\forall n=1,2,\cdots
\end{equation}
with initialisation condition
\begin{equation}
\Phi_0=\left(
\begin{array}{c}
\lambda_0\\
s_0
\end{array}
\right).
\end{equation}
It is straightforward to observe that for any $n$ the system (\ref{msystem}) has the solution
\begin{equation}
\Phi_n=T^n\Phi_{0}.
\end{equation}
Since for any $n$ the functions $\lambda_n$ and $s_n$ are real-analytic on the same interval $I$ where $s_0$ and $\lambda_0$ are real-analytic, we can expand the vector-valued function $\Phi_n$ into power series about some $x_0\in I$ with radius of convergence $R$ up to the nearest pole of $\lambda_0$ and $s_0$. Likewise, we can expand the matrix $A$ into a matrix power series about $x_0$ with radius of convergence $R$. More precisely, if we substitute
\begin{equation}
\Phi_n(x)=\sum_{m=0}^\infty\mathfrak{C}_{m,n}(x-x_0)^m,\quad
\mathfrak{A}(x)=\sum_{k=0}^\infty\mathfrak{A}_k(x-x_0)^k,\quad|x-x_0|<R
\end{equation}
with $\mathfrak{C}_{m,n}\in\mathbb{C}^2$ and $\mathfrak{A}_k\in\mathbb{C}^{2\times 2}$ into (\ref{msystem}) and we make use of the Cauchy product of power series together with a trivial shift of indices,  we end up with the following recurrence relation
\begin{equation}\label{rcof}
\mathfrak{C}_{m,n+1}=(m+1)\mathfrak{C}_{m+1,n}+\sum_{\ell=0}^m \mathfrak{A}_{m-\ell}\mathfrak{C}_{\ell,n}\quad\forall~m,n=0,1,\cdots.
\end{equation}
Given $x_0\in I$ and for each fixed $n$ we can recursively compute the coefficients $\mathfrak{C}_{m,n}$ for $m=1,2,\cdots$ from (\ref{rcof}). For instance, for $n=0$ equation (\ref{rcof}) allows to compute all coefficients in the Taylor expansion for $\Phi_1$ and so on. This shows the uniqueness of the solution $\Phi_n$ to (\ref{msystem}) for any $n$.

\section{The link between AIM and the continued fraction method}\label{Abt3}
In this section, we show how the unknown function $\alpha$ appearing in (\ref{Riccati}) can be represented in terms of a certain continued fraction. Furthermore, we also discuss the convergence problem of the latter. To do that, we follow a different approach {\bf{to}} the one outlined in Section~\ref{Intr}. If we differentiate (\ref{ode1}) one time with respect to $x$, we get
\begin{equation}\label{210}
y^{1+2}=\lambda_0 y^{''}+(\lambda_0^{'}+s_0)y^{'}+s_0^{'}y.
\end{equation}
On the other hand, we can solve (\ref{ode1}) with respect to $y$ to get
\begin{equation}
y=\frac{y^{''}}{s_0}-\frac{\lambda_0}{s_0}y^{'}
\end{equation}
which replaced into (\ref{210}) gives
\begin{equation}\label{321}
y^{(1+2)}=p_1 y^{''}+q_1 y^{'}
\end{equation}
with
\begin{equation}
p_1=\lambda_0+\frac{s_0^{'}}{s_0},\quad q_1=s_0+\lambda_0^{'}-\lambda_0\frac{s_0^{'}}{s_0}.
\end{equation}
Differentiating once more (\ref{321}) and substituting there
\begin{equation}
y^{'}=\frac{y^{(1+2)}}{q_1}-\frac{p_1}{q_1}y^{''}
\end{equation}
yields
\begin{equation}
y^{(2+2)}=p_2 y^{(1+2)}+q_2 y^{''}
\end{equation}
with
\begin{equation}
p_2=p_1+\frac{q_1^{'}}{q_1},\quad q_2=q_1+p_1^{'}-p_1\frac{q_1^{'}}{q_1}y^{''}.
\end{equation}
At this point, if we proceed iteratively, we obtain for all $n\geq 0$
\begin{equation}\label{box}
y^{(n+2)}=p_n(x)y^{(n+1)}+q_n(x)y^{(n)}
\end{equation}
where
\begin{eqnarray}
p_n&=&p_{n-1}+\frac{q^{'}_{n-1}}{q_{n-1}},\label{RICREL1}\\
q_n&=&q_{n-1}+p^{'}_{n-1}-p_{n-1}\frac{q^{'}_{n-1}}{q_{n-1}},\label{RICREL2}
\end{eqnarray}
with $p_0:=\lambda_0$, $q_0:=s_0$ and $p_{-1}=q_{-1}=0$. Equation (\ref{box}) is a second order linear ODE for the unknown function $y^{(n)}$ and, therefore, for each $n\geq 0$ we can associate to it a Riccati equation of the form
\begin{equation}\label{boxbox}
\varphi^{'}_n=\varphi_n^2(x)+p_n(x)\varphi_n(x)-q_n(x),\quad \varphi_n:=-\frac{y^{(n+1)}}{y^{(n)}}.
\end{equation}
It is straightforward to observe that in the case $n=0$ equation (\ref{boxbox}) coupled with the conditions $p_0:=\lambda_0$ and $q_0:=s_0$ reduces to the Riccati equation (\ref{Riccati}) and hence, we can make the identification $\varphi_0=\alpha$. If we rewrite (\ref{box}) as
\begin{equation}\label{cic}
\frac{y^{(n+2)}}{y^{(n+1)}}=p_n(x)+q_n(x)\frac{y^{(n)}}{y^{(n+1)}}
\end{equation}
and we take into account that $\varphi_{n+1}=-y^{(n+2)}/y^{(n+1)}$, then (\ref{cic}) becomes
\begin{equation}
\frac{q_n}{\varphi_n}=p_n+\varphi_{n+1}
\end{equation}
from which the following recurrence relation for $\varphi_n$ emerges, namely
\begin{equation}\label{threedots}
\varphi_n(x)=\frac{q_n(x)}{p_n(x)+\varphi_{n+1}(x)}\quad\forall n\geq 0.
\end{equation}
Finally, the solution of the Riccati equation (\ref{Riccati}) can be expressed as a continued fraction by means of (\ref{threedots}) as follows
\begin{equation}\label{CF}
\alpha(x)=\varphi_0(x)=\frac{q_0(x)}{p_0(x)+\varphi_{1}(x)}=\cfrac{q_0(x)}{p_0(x)+\cfrac{q_1(x)}{p_1(x)+\cfrac{q_2(x)}{p_2(x)+\dots}}}=\foo_{n=0}^\infty\left(\frac{q_n(x)}{p_n(x)}\right).
\end{equation}
At this point a couple or remarks are in order.
\begin{enumerate}
\item
If for some positive integer $N$ we have $q_N(x)=0$ for all $x\in I$, then (\ref{threedots}) induces a termination into the continued fraction (\ref{CF}) and by means of (\ref{threedots}) and (\ref{CF}) it follows that
\begin{equation}
\alpha(x)=\cfrac{q_0(x)}{p_0(x)+\cfrac{q_1(x)}{p_1(x)+\cfrac{q_2(x)}{\ddots p_{N-2}(x)+\cfrac{q_{N-1}(x)}{p_{N-1}(x)}}}}.
\end{equation}
\item
If no termination condition occurs, a careful analysis of the convergence/divergence property of the continuous fraction (\ref{CF}) is mandatory.
\end{enumerate}

\section{A series representation for the function $\alpha$}
In the following, the N-th approximant of (\ref{CF}) is defined as
\begin{equation}\label{mi}
C_N(x)=\frac{A_N(x)}{B_N(x)}=\foo_{n=0}^N\left(\frac{q_n(x)}{p_n(x)}\right),
\end{equation}
where the sequence of functions $\left(A_N(x)\right)_{n\in\mathbb{N}}$ and $\left(B_N(x)\right)_{n\in\mathbb{N}}$ are the $N$-th partial numerator and the N-th partial denominator of the continued fraction (\ref{CF}), respectively. We will always assume that $A_N/B_N$ is in reduced form, i.e. $A_N$ and $B_N$ have no common factors. Furthermore, we say that the continued fraction (\ref{CF}) converges pointwise to $\alpha$ on some interval $I$ of the real line if
\begin{equation}\label{limozzo}
\lim_{N\to\infty}C_N(x)=\alpha(x)\quad\forall~x\in I.
\end{equation}
By observing that
\begin{eqnarray}
C_0(x)&=&\frac{A_0(x)}{B_0(x)}=\frac{q_0(x)}{p_0(x)}=\frac{s_0(x)}{\lambda_0(x)}=\alpha_0(x),\\
C_1(x)&=&\frac{A_1(x)}{B_1(x)}=\frac{p_1(x)q_0(x)}{p_0(x)p_1(x)+q_1(x)}=\frac{s_0^{'}(x)+s_0(x)\lambda_0(x)}{\lambda_0^{'}(x)+\lambda_0^2(x)+s_0(x)}=\frac{s_1(x)}{\lambda_1(x)}=\alpha_1(x),\\
&\vdots&
\end{eqnarray}
we discover that the sequence of approximants of (\ref{CF}) is equivalent to the sequence $(\alpha_n)_{n\in\mathbb{N}}$ appearing in the AIM method. Hence, we proved that the AIM condition (\ref{ascond}) introduced by \cite{Ciftci} is equivalent to require that $\alpha$ is given by the $n$-th approximant of the continued fraction (\ref{CF}) for $n$ sufficiently large. On the other hand, if (\ref{CF}) or, equivalently, the sequence of its approximants converges pointwise on some interval $I$, then, the AIM asymptotic condition (\ref{equiv}) can be replaced by the condition
\begin{equation}\label{nic}
\lim_{n\to\infty}\left[C_n(x)-C_{n-1}(x)\right]=0\quad\forall~x\in I.
\end{equation}
In the case a termination behaviour occurs in a problem amenable to the treatment with the AIM method, the termination condition (\ref{term}) can be replaced by the equivalent condition 
\begin{equation}
q_N(x)=0\quad\forall~x\in I
\end{equation}
with $N$ some positive integer. The next theorem will allow us to reformulate the condition (\ref{nic}) as a certain limit involving only the terms of the sequence of the partial denominators of the continued fraction (\ref{CF}) and those belonging to sequence $(q_n)_{n\in\mathbb{N}}$ defined through the recurrence relation (\ref{RICREL1}) and (\ref{RICREL2}).
\begin{theorem}
Consider the continued fraction $\foo_{n=0}^N\left(q_n(x)/p_n(x)\right)$ whose $n$-th approximant is given by (\ref{mi}). Then, the $n$-th partial numerator $A_n(x)$ and denominator $B_n(x)$ of the continued fraction satisfy the three-term recurrence relation
\begin{eqnarray}
A_n(x)&=&p_n(x)A_{n-1}(x)+q_n(x)A_{n-2}(x),\quad A_{-1}(x)=0,~A_{-2}(x)=1\quad\forall x\in I,\label{kdo1}\\
B_n(x)&=&p_n(x)B_{n-1}(x)+q_n(x)B_{n-2}(x),\quad B_{-1}(x)=1,~B_{-2}(x)=0\quad\forall x\in I\label{kdo2}
\end{eqnarray}
for all $n\geq 0$.
\end{theorem}
\begin{proof}
We proceed by induction. First of all, we observe that for $n=0$ equations (\ref{kdo1}) and (\ref{kdo2}) reproduce correctly
\begin{eqnarray}
A_0(x)&=&p_0(x)A_{-1}(x)+q_0(x)A_{-2}(x)=q_0(x),\\
B_0(x)&=&p_0(x)B_{-1}(x)+q_0(x)B_{-2}(x)=p_0(x).
\end{eqnarray}
In the case $n=1$, we obtain as expected
\begin{eqnarray}
A_1(x)&=&p_1(x)A_{0}(x)+q_1(x)A_{-1}(x)=p_1(x)q_0(x),\\
B_1(x)&=&p_1(x)B_{0}(x)+q_1(x)B_{-1}(x)=p_0(x)p_1(x)+q_1(x).
\end{eqnarray}
Let us assume that (\ref{kdo1}) and (\ref{kdo2}) are true for $n=m$. Then, we have
\begin{eqnarray}
A_m(x)&=&p_m(x)A_{m-1}(x)+q_m(x)A_{m-2}(x),\label{drug}\\
B_m(x)&=&p_m(x)B_{m-1}(x)+q_m(x)B_{m-2}(x).\label{tret}
\end{eqnarray}
From
\begin{eqnarray}
C_{m+1}(x)&=&\foo_{n=0}^{m+1}\left(\frac{q_n(x)}{p_n(x)}\right)
=\cfrac{q_0(x)}{p_0(x)+\cfrac{q_1(x)}{p_1(x)+\cfrac{q_2(x)}{\ddots p_{m-1}(x)+\cfrac{q_{m}(x)}{p_{m}(x)+\frac{q_{m+1}(x)}{p_{m+1}(x)}}}}},\\
C_{m}(x)&=&\cfrac{q_0(x)}{p_0(x)+\cfrac{q_1(x)}{p_1(x)+\cfrac{q_2(x)}{\ddots p_{m-1}(x)+\frac{q_{m}(x)}{p_{m}(x)}}}}
\end{eqnarray}
we realize that $C_{m+1}(x)$ can be obtained from $C_m(x)$ by replacing there $p_m(x)$ with $p_m(x)+(q_{m+1}(x)/p_{m+1}(x))$. Hence, we can write
\begin{equation}
C_{m+1}(x)=\frac{A_{m+1}(x)}{B_{m+1}(x)}=\frac{\widehat{A}_{m}(x)}{\widehat{B}_{m}(x)},
\end{equation}
with
\begin{eqnarray}
\widehat{A}_{m}(x)&=&\frac{p_{m+1}(x)A_m(x)+q_{m+1}(x)A_{m-1}(x)}{p_{m+1}(x)},\label{aa}\\
\widehat{B}_{m}(x)&=&\frac{p_{m+1}(x)B_m(x)+q_{m+1}(x)B_{m-1}(x)}{p_{m+1}(x)},\label{bb}
\end{eqnarray}
where we made use of (\ref{drug}) and (\ref{tret}). Finally, by means of (\ref{aa}) and (\ref{bb}) we get
\begin{equation}
\frac{A_{m+1}(x)}{B_{m+1}(x)}=\frac{p_{m+1}(x)A_m(x)+q_{m+1}(x)A_{m-1}(x)}{p_{m+1}(x)B_m(x)+q_{m+1}(x)B_{m-1}(x)}
\end{equation}
and the proof is completed.~$\Box$
\end{proof}
A formula for the $N$-th approximant can be found multiplying (\ref{kdo1}) by $B_{n-1}(x)$, (\ref{kdo2}) by $A_{n-1}(x)$ followed by subtraction of both. This procedure leads to
\begin{equation}\label{ena}
A_n(x)B_{n-1}(x)-A_{n-1}(x)B_n(x)=-q_n(x)\left[A_{n-1}(x)B_{n-2}(x)-A_{n-2}(x)B_{n-1}(x)\right].
\end{equation}
Let $v_n(x)=A_n(x)B_{n-1}(x)-A_{n-1}(x)B_n(x)$. Using the initial conditions in (\ref{kdo1}) and (\ref{kdo2}) yields $v_{-1}(x)=-1$. Furthermore, (\ref{ena}) can be rewritten as
\begin{equation}\label{dva}
v_n(x)=-q_n(x)v_{n-1}(x)
\end{equation}
and its solution is
\begin{equation}\label{tri}
v_n(x)=(-1)^n\prod_{j=0}^{n}q_j(x)\quad\forall n\geq 0
\end{equation}
with $v_{-1}(x)=-1$. At this point, by means of (\ref{mi}) and the definition of $v_n$ we observe that
\begin{equation}
\lim_{n\to\infty}\left[C_n(x)-C_{n-1}(x)\right]=\lim_{n\to\infty}\frac{A_n(x)B_{n-1}(x)-A_{n-1}(x)B_n(x)}{B_n(x)B_{n-1}(x)}=
\lim_{n\to\infty}\frac{v_n(x)}{B_n(x)B_{n-1}(x)}.
\end{equation}
It is straightforward to verify with the help of (\ref{dva}) and (\ref{tri}) that the AIM condition (\ref{nic}) becomes
\begin{equation}\label{erg}
\lim_{n\to\infty}\frac{(-1)^n\prod_{j=0}^{n}q_j(x)}{B_n(x)B_{n-1}(x)}=0\quad\forall x\in I.
\end{equation}
The above formula clearly shows that the termination behaviour occurs when one of the $q_j$ vanishes. Furthermore, the question of whether or not the limit in (\ref{erg}) vanishes on the interval $I$ may be answered by constructing uniform asymptotic representations for $B_n(x)$ solution of the recurrence relation (\ref{tret}). If we divide (\ref{tri}) by $B_n(x)B_{n-1}(x)$ and introduce the difference operator $\Delta x_n:=x_{n+1}-x_n$ \cite{Elady} , we get
\begin{equation}\label{stiri}
\Delta C_{n-1}(x)=\frac{A_n(x)}{B_n(x)}-\frac{A_{n-1}(x)}{B_{n-1}(x)}=\frac{(-1)^n\prod_{j=0}^{n}q_j(x)}{B_n(x)B_{n-1}(x)}\quad\forall n\geq 0
\end{equation}
and applying the antidifference operator $\Delta^{-1}$ defined in $(2.1.16)$ in \cite{Elady} yields
\begin{equation}\label{enajst}
C_{n-1}(x)=\sum_{k=0}^{n-1}(-1)^k\frac{q_0(x)\cdots q_k(x)}{B_k(x)B_{k-1}(x)}\quad\forall n\geq 1.
\end{equation}
Combining (\ref{enajst}) with (\ref{limozzo}) and letting $n\to\infty$ leads to the following formal representation of the continued fraction (\ref{CF}) or equivalently for the unknown function $\alpha$, namely
\begin{equation}\label{zve}
\alpha(x)=\foo_{n=0}^\infty\left(\frac{q_n(x)}{p_n(x)}\right)=\sum_{n=0}^\infty(-1)^n\frac{q_0(x)\cdots q_n(x)}{B_n(x)B_{n-1}(x)}.
\end{equation}
In order to formulate a necessary and sufficient condition for the continued fraction (\ref{zve}) to converge pointwise to some function on a given interval of the real line, we first recall that two continued fractions are said to be equivalent (denoted by the symbol $\approx$) if both have the same sequence of approximants.

\begin{theorem}\label{Ena}
For all $n\in\mathbb{N}\cup\{0\}$ and for all $x\in\mathfrak{I}\subseteq\mathbb{R}$ suppose that $p_n(x)>0$ and $q_n(x)$ does not vanish. Further assume that $\mathcal{I}=\mathfrak{I}\cap I\neq\emptyset$ where $I$ is the interval introduced in (\ref{ode1}). Then, the continued fraction (\ref{zve}) converges pointwise on $\mathcal{I}$ to some function $\alpha$ if and only if the series 
\begin{equation}
\sum_{n=0}^\infty p_n(x)
\end{equation}
diverges for every $x\in\mathcal{I}$.
\end{theorem}
\begin{proof}
Fix any $x_0\in\mathcal{I}$. We start with the observation that the equivalence relation
\begin{equation}
\foo_{n=0}^\infty\left(\frac{q_n(x_0)}{p_n(x_0)}\right)\approx\foo_{n=0}^\infty\left(\frac{\gamma_{n-1}\gamma_n q_n(x_0)}{\gamma_n p_n(x_0)}\right)
\end{equation}
holds for any sequence of non zero complex numbers $\gamma_{-1}=1$, $\gamma_0$, $\gamma_1\cdots$ \cite{Elady}. In particular, if we choose the sequence $(\gamma_n)_{n\in\mathbb{N}}$ such that $\gamma_{n-1}=1/q_n$, then
\begin{equation}
\foo_{n=0}^\infty\left(\frac{q_n(x_0)}{p_n(x_0)}\right)\approx q_0(x_0)\foo_{n=0}^\infty\left(\frac{1}{p_n(x_0)}\right).
\end{equation}
From (\ref{zve}) it follows that
\begin{equation}\label{sum}
\foo_{n=0}^\infty\left(\frac{1}{p_n(x_0)}\right)=\sum_{n=0}^\infty S_n,\quad S_n:=\frac{(-1)^n}{B_{n-1}(x_0)B_n(x_0)}
\end{equation}
with $B_n(x_0)$ given by the recurrence relation (\ref{tret}) with $q_n=1$ for all $n$, namely
\begin{equation}\label{tilde}
B_n(x_0)=p_n(x_0)B_{n-1}(x_0)+B_{n-2}(x_0),\quad B_{-2}(x_0)=0,\quad B_{-1}(x_0)=1.
\end{equation}
Note that $p_n(x)>$ for all $n\in\mathbb{N}\cup\{0\}$ and for all $x\in\mathcal{I}$ ensures that all terms $B_n(x_0)$ in the above recurrence relation are positive. From (\ref{sum}) we realize that (\ref{zve}) converges whenever the alternating series in (\ref{sum}) converges. To this purpose, we observe that
\begin{equation}
B_{n+1}(x_0)=p_{n+1}(x_0)B_n(x_0)+B_{n-1}(x_0)>B_{n-1}(x_0).
\end{equation}
This implies that $B_n(x_0)B_{n+1}(x_0)>B_n(x_0)B_{n-1}(x_0)$ and hence,
\begin{equation}
\frac{1}{B_n(x_0)B_{n+1}(x_0)}<\frac{1}{B_n(x_0)B_{n-1}(x_0)}
\end{equation}
for all $n\in\mathbb{N}\cup\{0\}$. Thus, the sequence $(|S_n|)$ is monotonically decreasing. Therefore, the series in (\ref{sum}) converges provided that the following limit diverges
\begin{equation}\label{ahah}
\lim_{n\to\infty}B_{n-1}(x_0)B_n(x_0).
\end{equation}
Furthermore, from (\ref{tilde}) we have
\begin{eqnarray}
B_1(x_0)&=&1+p_0(x_0)p_1(x_0)>1,\\
B_2(x_0)&=&p_0(x_0)+p_2(x_0)+p_0(x_0)p_1(x_0)p_2(x_0)>p_0(x_0),\\
B_3(x_0)&=&1+p_0(x_0)p_1(x_0)+p_0(x_0)p_3(x_0)+p_2(x_0)p_3(x_0)+p_0(x_0)p_1(x_0)p_2(x_0)p_3(x_0)>1,\\
&\vdots&\nonumber
\end{eqnarray}
Hence, $B_n(x_0)>\min{\{1,p_0\}}=:\mu$ and we conclude that
\begin{equation}\label{trideset}
B_{n-1}(x_0)B_n(x_0)>\mu^2\sum_{k=0}^n p_k(x_0).
\end{equation}
This shows that if $\sum_{n=0}^\infty p_n(x_0)$ diverges, the limit in (\ref{ahah}) diverges as well and the continued fraction converges. Suppose now that $\sum_{n=0}^\infty p_n(x_0)$ converges. Then, (\ref{tilde}) implies that
\begin{equation}
B_{n-1}(x_0)+B_n(x_0)=B_{n-2}(x_0)+\left[1+p_n(x_0)\right]B_{n-1}(x_0)<\left[1+p_n(x_0)\right]\left[B_{n-1}(x_0)+B_{n-2}(x_0)\right]
\end{equation}
and by induction
\begin{equation}
B_{n-1}(x_0)+B_{n}(x_0)<\prod_{k=0}^n\left[1+p_k(x_0)\right]<\prod_{k=0}^n e^{p_k(x_0)}<\mbox{exp}\left(\sum_{n=0}^\infty p_n(x_0)\right).
\end{equation}
Since $\sum_{n=0}^\infty p_n(x_0)$ converges, let $L$ denote the corresponding limit. This implies that
\begin{equation}
B_{n-1}(x_0)+B_{n}(x_0)<e^L
\end{equation}
and hence,
\begin{equation}
B_{n-1}(x_0)B_{n}(x_0)<B^2_{n-1}(x_0)+B_{n}^2(x_0)+2B_{n-1}(x_0)B_{n}(x_0)=
\left[B_{n-1}(x_0)+B_{n}(x_0)\right]^2<e^{2L}.
\end{equation}
This results violates the divergence condition for (\ref{ahah}) and (\ref{zve}) must diverge at $x_0$. Since the choice of $x_0\in\mathcal{I}$ was arbitrary, the result is valid for any $x\in\mathcal{I}$.~$\Box$
\end{proof}
In order to probe that the function $\alpha$ can be interpreted as a certain term of a minimal solution of a recurrence relation we derive hereafter, we briefly recall that a solution $\psi_n$ of the difference equation
\begin{equation}
\mathfrak{x}_{n+2}+\mathfrak{f}_1(n)\mathfrak{x}_{n+1}+\mathfrak{f}_2(n)\mathfrak{x}_{n}=0
\end{equation}
is called minimal or subdominant if
\begin{equation}
\lim_{n\to\infty}\frac{\psi_n}{\mathfrak{x_n}}=0
\end{equation}
for any solution $\mathfrak{x}_n$ that is not a multiple of $\psi_n$.
\begin{theorem}\label{TH3}
Let $\mathcal{I}$ be defined as Theorem~IV.2. Consider the recurrence relation
\begin{equation}\label{3rr}
\mathfrak{x}_{n}(x)-p_n(x)\mathfrak{x}_{n-1}(x)-q_n(x)\mathfrak{x}_{n-2}(x)=0
\end{equation}
with $q_n(x)\neq 0$ for all $x\in\mathcal{I}$ and for every $n\in\mathbb{N}\cup\{-1,0\}$. If for any $x\in\mathcal{I}$ (\ref{3rr}) has a minimal solution $\psi_n(x)$ with $\psi_{-1}(x)\neq 0$ for all $x\in\mathcal{I}$, then the continued fraction (\ref{zve}) converges pointwise to $\alpha$ on the interval $\mathcal{I}$.
\end{theorem}
\begin{proof}
Suppose that (\ref{3rr}) has a minimal solution $\psi_n(x)$ with $\psi_{-1}(x)\neq 0$ for all $x\in\mathcal{I}$. If we divide (\ref{3rr}) by $\mathfrak{x}_{n-1}$ and introduce $\mathfrak{y}_n(x):=\mathfrak{x}_{n}(x)/\mathfrak{x}_{n-1}(x)$ with $\mathfrak{x}_{n-1}(x)\neq 0$ for all $x\in\mathcal{I}$, we obtain 
\begin{equation}
\mathfrak{y}_n(x)-p_n(x)=\frac{q_n(x)}{\mathfrak{y}_{n-1}(x)}
\end{equation}
and hence,
\begin{equation}
\mathfrak{y}_{n-1}(x)=\frac{q_n(x)}{-p_n(x)+\mathfrak{y}_n(x)}.
\end{equation}
A recursive application of the formula above gives
\begin{equation}
\mathfrak{y}_{n-1}(x)=\cfrac{q_n(x)}{-p_n(x)+\cfrac{q_{n+1}(x)}{-p_{n+1}(x)+\cfrac{q_{n+2}(x)}{ -p_{n+2}(x)+\cdots}}}.
\end{equation}
Taking into account that $\mathfrak{y}_{n-1}(x):=\mathfrak{x}_{n-1}(x)/\mathfrak{x}_{n-2}(x)$, the case $n=0$ yields
\begin{equation}\label{bubu}
\mathfrak{y}_{-1}(x)=\frac{\mathfrak{x}_{-1}(x)}{\mathfrak{x}_{-2}(x)}=
\cfrac{q_0(x)}{-p_0(x)+\cfrac{q_{1}(x)}{-p_{1}(x)+\cfrac{q_{2}(x)}{ -p_{2}(x)+\cdots}}}=
\foo_{n=0}^\infty\left(\frac{q_n(x)}{-p_n(x)}\right)
\end{equation}
with $n$-th approximant
\begin{equation}
\widehat{C}_n(x)=\frac{\widehat{A}_n(x)}{\widehat{B}_n(x)}.
\end{equation}
On the other hand, $\widehat{A}_n(x)$ and $\widehat{B}_n(x)$ are two linearly independent solutions of (\ref{3rr}) with initial conditions $\widehat{A}_{-1}(x)=0$, $\widehat{A}_{-2}(x)=1$ and $\widehat{B}_{-1}(x)=1$, $\widehat{B}_{-2}(x)=0$ for all $x\in\mathcal{I}$. Hence, without loss of generality let for some $\ell(x)\neq 0$ and for all $x\in\mathcal{I}$
\begin{equation}
\psi_n(x)=\widehat{A}_{n}(x)+\ell(x)\widehat{B}_{n}(x)\quad\forall n\geq-1.
\end{equation}
Note that for all $x\in\mathcal{I}$
\begin{equation}
\psi_{-1}(x)=\ell(x),\quad \psi_{-2}(x)=1.
\end{equation}
Furthermore, we have
\begin{equation}\label{ru}
\frac{\psi_n(x)}{\widehat{B}_{n}(x)}=\frac{\widehat{A}_n(x)}{\widehat{B}_n(x)}+\ell(x).
\end{equation}
Since $\psi_n(x)$ is a minimal solution of (\ref{3rr}) for each $x\in\mathcal{I}$, it follows that
\begin{equation}
\lim_{n\to\infty}\frac{\psi_n(x)}{\widehat{B}_{n}(x)}=0
\end{equation}
pointwise on $\mathcal{I}$. If we take $n\to\infty$ in (\ref{ru}), we obtain again pointwise on $\mathcal{I}$
\begin{equation}\label{elle}
\lim_{n\to\infty}\widehat{C}_n(x)=\lim_{n\to\infty}\frac{\widehat{A}_n(x)}{\widehat{B}_n(x)}=-\ell(x).
\end{equation}
On the other hand, from (\ref{bubu}) we have
\begin{equation}
\mathfrak{y}_{-1}(x)=\frac{\psi_{-1}(x)}{\psi_{-2}(x)}=\ell(x)
\end{equation}
which combined with (\ref{elle}) and $\psi_{-2}(x)=1$ gives
\begin{equation}
\lim_{n\to\infty}\widehat{C}_n(x)=-\psi_{-1}(x).
\end{equation}
Taking into account the general property that if a convergent continued fraction $\foo_{n=0}^\infty\left(\frac{a_n}{-b_n}\right)$ has limit $L$, then $\foo_{n=0}^\infty\left(\frac{a_n}{b_n}\right)$ has also limit $L$, it follows that
\begin{equation}
\lim_{n\to\infty}C_n(x)=-\lim_{n\to\infty}\widehat{C}_n(x)=\psi_{-1}(x)
\end{equation}
and the proof is completed.~$\Box$
\end{proof}
This result tells us that if the conditions in the theorem above are met, the sequence with terms $\alpha_n=s_n/\lambda_n$ in the AIM method will converge pointwise on $\mathcal{I}$ to some function $\alpha$ which can be interpreted as the term $\psi_{-1}(x)$ of the minimal solution $\psi_n(x)$ of the recurrence relation (\ref{3rr}). We investigate now some conditions ensuring the existence of a minimal solution. To this purpose, it is convenient to introduce the change of variable
\begin{equation}\label{KK}
\mathfrak{x}_n(x)=\left(-\frac{1}{2}\right)^{n-1}\prod_{j=n_0}^{n-2}p_j(x)\mathfrak{y}_n(x),
\end{equation}
which applied to (\ref{3rr}) leads to the recurrence relation
\begin{equation}\label{konec}
\mathfrak{y}_n(x)-2\mathfrak{y}_{n+1}(x)+\frac{4q_n(x)}{p_{n-1}(x)p_{n}(x)}\mathfrak{y}_{n+2}(x)=0.
\end{equation}
Suppose that the following limit exists pointwise everywhere on $\mathcal{I}_0\subseteq \mathcal{I}$, namely
\begin{equation}
q(x)=\lim_{n\to\infty}\frac{4q_n(x)}{p_{n-1}(x)p_{n}(x)}
\end{equation}
and define for every $x\in\mathcal{I}_0$
\begin{equation}\label{adef}
a_n(x)=\frac{4q_n(x)}{p_{n-1}(x)p_{n}(x)}-q(x).
\end{equation}
Then, the characteristic roots of (\ref{konec}) are given by
\begin{equation}
r_{\pm}(x)=1\pm\sqrt{1-q(x)}\quad\forall x\in\mathcal{I}_0.
\end{equation}
According to \cite{Elady} we have the following cases
\begin{enumerate}
\item
If $q(x)<1$ for all $x\in\mathcal{I}_1\subseteq\mathcal{I}$ and $p_n(x)$, $q_n(x)$ are real-valued functions, then $r_\pm(x)$ are real and distinct. Furthermore, we have the following subcases
\begin{enumerate}
\item
If $a_n(x)\to 0$ as $n\to\infty$ pointwise for any $x\in\mathcal{I}_2\subseteq\mathcal{I}$, then the Poincar$\acute{\mbox{e}}$ - Perron theorem ensures the existence of two linearly independent solutions $\mathfrak{y}_{\pm,n}(x)$ of (\ref{konec}) such that
\begin{equation}\label{alabarda}
\lim_{n\to\infty}\frac{\mathfrak{y}_{\pm,n+1}(x)}{\mathfrak{y}_{\pm,n}(x)}=r_{\pm}(x)\quad\forall x\in \mathcal{I}_0\cap\mathcal{I}_1\cap\mathcal{I}_2.
\end{equation}
Let $\mathcal{I}_3\subseteq\mathcal{I}$ be some interval of the real line where $|r_{-}(x)|<|r_+(x)|$ for any $x\in\mathcal{I}_3$, then we can construct some positive functions $f_{\pm}(x)$ such that
\begin{equation}
|r_{-}(x)|<f_{-}(x)<f_{+}(x)<|r_{+}(x)|\quad\forall x\in\mathcal{I}_3.
\end{equation}
Moreover, (\ref{alabarda}) implies that
\begin{equation}
\left|\frac{\mathfrak{y}_{-,n+1}(x)}{\mathfrak{y}_{-,n}(x)}\right|\leq f_{-}(x),\quad
\left|\frac{\mathfrak{y}_{+,n+1}(x)}{\mathfrak{y}_{+,n}(x)}\right|\geq f_{+}(x)\quad \forall x\in\bigcap_{k=0}^3\mathcal{I}_k
\end{equation}
from which
\begin{equation}
|\mathfrak{y}_{-,n}(x)|\leq f_{-}^n(x)|\mathfrak{y}_{-,0}(x)|,\quad
|\mathfrak{y}_{+,n}(x)|\geq f_{+}^n(x)|\mathfrak{y}_{+,0}(x)|
\end{equation}
and hence,
\begin{equation}
\lim_{n\to\infty}\left|\frac{\mathfrak{y}_{-,n}(x)}{\mathfrak{y}_{+,n}(x)}\right|=\left|\frac{\mathfrak{y}_{-,0}(x)}{\mathfrak{y}_{+,0}(x)}\right|
\lim_{n\to\infty}\left(\frac{f_-(x)}{f_+(x)}\right)^n=0.
\end{equation}
This shows that $\mathfrak{y}_{-,n}(x)$ is a minimal solution to (\ref{konec}) and by Theorem~\ref{TH3} the continued fraction (\ref{zve}) converges pointwise to $\alpha$ on the interval $\bigcap_{k=0}^3\mathcal{I}_k$.
\item
Suppose there exists some interval $\mathcal{I}_4\subseteq\mathcal{I}$ such that for any $x\in\mathcal{I}_4$
\begin{equation}
\sum_{n=n_0}^\infty |a_n(x)|<\infty.
\end{equation}
Then, Corollary 8.27 in \cite{Elady} allows to construct two linearly independent solutions of (\ref{konec}) having asymptotic behaviour
\begin{equation}
\mathfrak{y}_{\pm,n}(x)=r_{\pm}^n(x)\left[1+o(1)\right]
\end{equation}
and by means of (\ref{KK}) the recurrence relation (\ref{3rr}) has solutions
\begin{equation}
\mathfrak{x}_{\pm,n}(x)=\left(-\frac{1}{2}\right)^{n-1}\prod_{j=n_0}^{n-2}p_j(x)r_{\pm}^n(x)\left[1+o(1)\right]\quad
\forall x\in\bigcap_{k=0}^4\mathcal{I}_k.
\end{equation}
Taking into account that the inequality $|r_{-}(x)|<|r_+(x)|$ holds on $\mathcal{I}_3$, then
\begin{equation}
\lim_{n\to\infty}\left|\frac{\mathfrak{x}_{-,n}(x)}{\mathfrak{x}_{+,n}(x)}\right|=
\lim_{n\to\infty}\left|\frac{r_{-}(x)}{r_{+}(x)}\right|^n\left[1+o(1)\right]\quad\forall x\in\bigcap_{k=0}^4\mathcal{I}_k
\end{equation}
and we conclude that $\mathfrak{x}_{-,n}(x)$ is a minimal solution.
\item
If instead
\begin{equation}
\sum_{n=n_0}^\infty\left|a_n(x)\right|^2<\infty
\end{equation}
on an appropriate interval $\widehat{I}$, \cite{El} showed that (\ref{3rr}) has two linearly independent solutions
\begin{equation}
\mathfrak{x}_{\pm,n}(x)=\left(-\frac{1}{2}\right)^{n}\prod_{j=n_0}^{n-2}p_j(x)\prod_{j=n_0}^{n-1}\left[1\mp\sqrt{1-q(x)}\pm\frac{a_j(x)}{2\sqrt{1-q(x)}}\right]\left[1+o(1)\right].
\end{equation}
Then, $\mathfrak{x}_{-,n}(x)$ is a minimal solution and Theorem~\ref{TH3} can be used to conclude that the continued fraction (\ref{zve}) converges pointwise to $\alpha$.
\end{enumerate}
\item
This case is relevant when we need to consider those values of $x$ where $q(x)=1$ and hence, $r_\pm(x)=1$. As usual, let $x\in\mathcal{I}$ and $q(x)=1$ for all $x\in\widehat{\mathcal{I}_1}\subseteq\mathcal{I}$. Further suppose that
\begin{equation}
\sum_{n=n_0}^\infty n\left|a_n(x)\right|<\infty\quad\forall x\in\widehat{\mathcal{I}}_2\subseteq\mathcal{I}.
\end{equation}
Then, Coffman's Theorem (see Theorem 8.29 in \cite{Elady}) implies there exists two asymptotic solutions $\mathfrak{y}_{-,n}(x)\approx 1$ and $\mathfrak{y}_{+,n}(x)\approx n$, and hence, we can construct again a minimal solution.
\item
In the case $q(x)>1$ on an appropriate subset of $\mathcal{I}$, then $r_\pm(x)$ are complex conjugate. If in addition we require that
\begin{equation}
\sum_{n=n_0}^\infty\left|a_{n+1}(x)-a_n(x)\right|<\infty
\end{equation}
on some subset of $\mathcal{I}$, then, we can construct two linearly independent solutions as in \cite{El1}, namely
\begin{equation}
\mathfrak{y}_{\pm,n}(x)=\left[1+o(1)\right]\prod_{m=n_0}^{n-1}\beta_{\pm,m}(x),\quad \beta_{\pm,m}(x)=1\pm\sqrt{1-q(x)+a_n(x)}
\end{equation}
provided that
\begin{equation}
\Re{\sqrt{1-q(x)+a_n(x)}}\geq 0
\end{equation}
for a fixed branch of the square root, for instance $0\leq\arg{\sqrt{z}}<\pi$. Also in this case we can extract a minimal solution.
\item
Suppose that the coefficients in (\ref{3rr}) have the asymptotic expansions
\begin{equation}
p_n(x)\approx a(x)n^\sigma(x),\quad q_n(x)\approx b(x)n^{\tau(x)}
\end{equation}
with $a(x)b(x)\neq 0$ on some interval and $\sigma(x),\tau(x)$ real-valued functions. Then, we have the following criteria to check the existence of a minimal solution \cite{Kr,Per}
\begin{enumerate}
\item
If $\sigma(x)>\tau(x)/2$ on a subset of $\mathcal{I}$, we can find a fundamental set of solutions $\mathfrak{x}_{\pm,n}(x)$ such that the following limits hold pointwise
\begin{equation}
\lim_{n\to\infty}\frac{\mathfrak{x}_{-,n+1}(x)}{\mathfrak{x}_{-,n}(x)}=-a(x)n^{\sigma(x)},\quad
\lim_{n\to\infty}\frac{\mathfrak{x}_{+,n+1}(x)}{\mathfrak{x}_{+,n}(x)}=-\frac{b(x)}{a(x)}n^{\sigma(x)-\tau(x)}.
\end{equation}
Furthermore, $\mathfrak{x}_{+,n}(x)$ is a minimal solution.
\item
If at some point $\sigma(x)=\tau(x)/2$, let $r_\pm(x)$ denote the roots of the equation $\lambda^2+\sigma(x)\lambda+\tau(x)=0$ with $|r_+(x)|\geq|r_{-}(x)|$ on some subset of $\mathcal{I}$. Then, for $|r_+(x)|\neq|r_{-}(x)|$ the recurrence relation (\ref{3rr}) has a fundamental set of solutions $\mathfrak{x}_{\pm,n}(x)$ such that the following limits hold pointwise
\begin{equation}
\lim_{n\to\infty}\frac{\mathfrak{x}_{-,n+1}(x)}{\mathfrak{x}_{-,n}(x)}=r_{-}(x)n^{\sigma(x)},\quad
\lim_{n\to\infty}\frac{\mathfrak{x}_{+,n+1}(x)}{\mathfrak{x}_{+,n}(x)}=r_{+}(x)n^{\tau(x)}.
\end{equation}
Furthermore, $\mathfrak{x}_{+,n}(x)$ is a minimal solution.
\end{enumerate}
\item
Suppose that the coefficients in (\ref{3rr}) have asymptotic expansions
\begin{equation}
p_n(x)\approx\sum_{j=0}^\infty\frac{a_j(x)}{n^j},\quad
q_n(x)\approx\sum_{j=0}^\infty\frac{b_j(x)}{n^j}
\end{equation}
with $b_0(x)\neq$ on some interval. The characteristic equation associated to (\ref{3rr}) is $r^2+a_0(x)r+b_0(x)=0$ with roots
\begin{equation}
r_\pm(x)=-\frac{a_0(x)}{2}\pm\sqrt{\frac{a_0^2(x)}{4}-b_0(x)}.
\end{equation}
At this point, the Birkhoff-Adams theorem \cite{BA1,BA2,BA3,BA4,LI,Wong} can be used in order to construct asymptotic expansions of the fundamental set of solutions to (\ref{3rr}) allowing to determine whether or not a minimal solution exists. We have the following cases
\begin{enumerate}
\item
If $a_0^2(x)\neq 4b_0(x)$ on some subset of $\mathcal{I}$, there are two linearly independent solutions $\mathfrak{x}_{\pm,n}(x)$ given by
\begin{equation}
\mathfrak{x}_{+,n}(x)=r_\pm^n(x)n^{\alpha_\pm(x)}\sum_{k=0}^\infty\frac{c_{\pm,k}(x)}{n^k},\quad\alpha_\pm(x)=\frac{a_1(x)r_{\pm(x)}+b_1(x)}{a_0(x)r_\pm(x)+2b_0(x)}
\end{equation}
for those values of $x$ such that the denominator in the expression for $\alpha_\pm(x)$ does not vanish. Furthermore, the coefficients $c_{\pm,k}$ satisfy the recurrence relation
\begin{equation}
\sum_{j=0}^{s-1}\left[r_\pm^2(x)2^{s-j}\left(
\begin{array}{c}
\alpha_\pm(x)-j\\
s-j
\end{array}
\right)+r_{\pm}(x)\sum_{k=j}^s\left(
\begin{array}{c}
\alpha_\pm(x)-j\\
k-j
\end{array}
\right)a_{s-k}(x)+b_{s-j}(x)\right]c_{\pm,j}(x)=0,\quad c_{\pm,0}(x)=1.
\end{equation}
\item
At those points where $r_{-}(x)=r_{+}(x)=r(x)$ with $r(x)=-a_0(x)/2$ not a root of $a_1(x)r(x)+b_1(x)=0$, (\ref{3rr}) has two linearly independent solutions of the form
\begin{equation}
\mathfrak{x}_{+,n}(x)=r^n(x)e^{\gamma_\pm(x)\sqrt{n}}n^{\widetilde{\alpha}_\pm(x)}\sum_{j=0}^\infty\frac{c_{\pm,j}(x)}{n^{j/2}}
\end{equation}
with
\begin{eqnarray}
\widetilde{\alpha}(x)&=&\frac{1}{4}+\frac{b_1(x)}{2b_0(x)},\quad\gamma_\pm(x)=\pm 2\sqrt{\frac{a_0(x)a_1(x)-2b_1(x)}{2b_0(x)}},\quad c_{\pm,0}=1,\\
c_{\pm,1}(x)&=&\frac{1}{24b_0^2(x)\gamma_\pm(x)}\left(
a_0^2(x)a_1^2(x)-24a_0(x)a_1(x)b_0(x)+8a_0(x)a_1(x)b_1(x)\right.\\
&&\left.-24a_0(x)a_2(x)b_0(x)-9b_0^2(x)-32b_1^2(x)+24b_0(x)b_1(x)+48b_0(x)b_2(x)
\right)
\end{eqnarray}
\item
At those points where $r_{-}(x)=r_{+}(x)=r(x)$ and $2b_1(x)=a_0(x)a_1(x)$, we consider instead the equation
\begin{equation}\label{ekv}
\widetilde{\alpha}(\widetilde{\alpha}-1)r^2(x)+\left[a_1(x)\widetilde{\alpha}+a_2(x)\right]r(x)+b_2(x)=0.
\end{equation}
By $\widetilde{\alpha}_\pm$ with $\Re{\widetilde{\alpha}_+}\geq\Re{\widetilde{\alpha}_-}$ we denote the roots of (\ref{ekv}). Then, we have the following three scenarios
\begin{enumerate}
\item
For those values of $x$ such that $\widetilde{\alpha}_+-\widetilde{\alpha}_-\neq 0,1,2,\cdots$ the equation (\ref{3rr}) has two linearly independent solutions of the form
\begin{equation}\label{zak}
\mathfrak{x}_{\pm,n}(x)=r^n(x)n^{\widetilde{\alpha}_\pm(x)}\sum_{j=0}^\infty\frac{c_{\pm,j}(x)}{n^{j}}.
\end{equation}
\item
For those values of $x$ such that $\widetilde{\alpha}_+-\widetilde{\alpha}_-= 1,2,\cdots$, then the equation (\ref{3rr}) has one solution $\mathfrak{x}_{-,n}(x)$ given by (\ref{zak}) and another solution
\begin{equation}
\mathfrak{x}_{+,n}(x)=z_n(x)+c\mathfrak{x}_{-,n}(x)\ln{n},
\end{equation}
where $c$ is a constant that may be zero and
\begin{equation}\label{zak1}
z_n(x)=r^n(x)n^{\widetilde{\alpha}_+(x)}+\sum_{s=0}^\infty\frac{d_s}{n^s}.
\end{equation}
\item
For those values of $x$ such that $\widetilde{\alpha}_+=\widetilde{\alpha}_-$, then the equation (\ref{3rr}) has one solution $\mathfrak{x}_{-,n}(x)$ given by (\ref{zak}) and a second solution represented by (\ref{zak1}) with
\begin{equation}
z_n(x)=r^n(x)n^{\widetilde{\alpha}_-(x)-k_0+2}+\sum_{s=0}^\infty\frac{d_s}{n^s}\sum_{s=0}^\infty\frac{d_s}{n^s},
\end{equation}
where $k_0$ is a positive integer such that $k_0\geq 3$.
\end{enumerate}
\end{enumerate}
\end{enumerate}
We observe that the above conditions ensuring the existence of a minimal solution may, in general, have the effect to shrink the interval over which $x$ must be chosen when we use the Asymptotic Iterative Method . The next result narrows the interval over which $x$ must be taken so that the continued fraction (\ref{sum}) and hence, (\ref{zve}) converge pointwise.
\begin{theorem}
Let $\mathfrak{K}:=\{x\in I~|~p_n(x)\geq 1\}$. Then,the AIM condition (\ref{nic}) is satisfied for all $x\in\mathfrak{K}$.
\end{theorem}
\begin{proof}
Since (\ref{sum}) and (\ref{zve}) are equivalent, there is no loss in generality in considering the sequence of approximants of (\ref{sum}). By (\ref{enajst}) we have for the $n$-th and $(n-1)$-th approximants
\begin{equation}
C_n(x)-C_{n-1}(x)=\frac{(-1)^n}{B_n(x)B_{n-1}(x)}
\end{equation}
with $B_n(x)$ satisfying the recurrence relation in (\ref{kdo2}). From (\ref{trideset}) we know that
\begin{equation}
B_n(x)B_{n-1}(x)>\mu^2\sum_{k=0}^n p_k(x),\quad\mu=\min_{x\in I}{\{1,p_0(x)\}}
\end{equation}
since $p_n(x)\geq 1$ for all $n$ ensures that $p_n(x)>0$ for all $n$. For $x\in\mathfrak{K}$ we have
\begin{equation}
\sum_{k=0}^n p_k(x)\geq n+1.
\end{equation}
Moreover, since all $p_k(x)$ are positive, we also have
\begin{equation}
\left|\sum_{k=0}^n p_k(x)\right|\geq n+1
\end{equation}
and therefore,
\begin{equation}
\left|B_n(x)B_{n-1}(x)\right|>(n+1)\mu^2
\end{equation}
from which it follows that
\begin{equation}
\left|C_n(x)-C_{n-1}(x)\right|<\frac{1}{(n+1)\mu^2}.
\end{equation}
Clearly, the above expression vanishes as $n\to\infty$ and the proof is completed.~$\Box$
\end{proof}

\begin{acknowledgements}
We thank the anonymous referees for the critical reading and useful comments which helped improve the quality of the paper.
\end{acknowledgements}


\begin{thebibliography}{999}
\bibitem{Ciftci}
H. Ciftci, R. L. Hall and N. Saad, {\it{Asymptotic iteration method for eigenvalue problems}}, J. Phys. A: Math. Gen. {\bf{36}}, 11807 (2003).
\bibitem{Ismail}
M. E. H. Ismail and N. Saad, {\it{The Asymptotic Iteration Method Revisited}}, J. Math. Phys. {\bf{61}}, 033501 (2020).
\bibitem{Cho1}
H. T. Cho, A.S. Cornell, J. Doukas and W. Naylor, {\it{Black hole quasinormal modes using the asymptotic iteration method}}, Class. Quantum Grav. {\bf{27}}, 155004 (2010).
\bibitem{Cho2}
H. T. Cho, A.S. Cornell, J. Doukas, T.-R. Huang and W. Naylor, {\it{A New Approach to Black Hole Quasinormal Modes: A Review of the Asymptotic Iteration Method}}, Adv. Math. Phys. {\bf{2012}} (2012) 281705.
\bibitem{Cama}
C. Camacho and H. Movasati, {\it{Remarks on a theorem of Perron}}, J. Differ. Equ. {\bf{260}}, 1465 (2016).
\bibitem{Matamala}
A.R. Matamala, F.A. Gutierrez and  J. Diaz-Vald$\acute{\mbox{e}}$s, {\it{A connection between the asymptotic iteration method and the continued fractions formalism}}, Phys. Lett. A{\bf{361}}, 16 (2007).
\bibitem{34}
D. Batic, M. Nowakowski and K. Redway, {\it{Some exact quasi-normal frequencies of a massless scalar field in the Schwarzschild space-time}}, Phys. Rev. D{\bf{98}}, 024017 (2018).
\bibitem{35}
R. Panosso Macedo, {\it{Comment on ”Some exact quasinormal frequencies of a massless scalar field in Schwarzschild space-time}}, Phys. Rev. D{\bf{99}}, 088501 (2019)
\bibitem{36}
D. Batic, M. Nowakowski and K. Redway, {\it{Reply to ”Comment on ‘Some exact quasinormal frequencies of a massless scalar field in Schwarzschild spacetime’ ”}}, Phys. Rev. D{\bf{99}}, 088502 (2019).
\bibitem{37}
H.Ciftci, R. L. Hall and N. Saad, {\it{Perturbation Theory in a framework of iteration methods}}, Phys. Lett. A{\bf{340}}, 388 (2005).
\bibitem{Elady}
S. Elaydi, {\it{An Introduction to Difference Equations}}, Springer Verlag (2005).
\bibitem{El}
S. Elaydi, {\it{Asymptotics for linear difference equations I: Basic theory}}, J. Difference Eqn. Appl. {\bf{5}}, 563 (1999).
\bibitem{El1}
S. Elaydi, {\it{Asymptotics for linear difference equations II: Application}} in {\it{New Trends in Difference Equations}}, Taylor \& Francis, London (2002).
\bibitem{Kr}
P. Kreuser, {\it{\"Uber das Verhalten der Integrale homogener linearer Differenzengleichungen im Unendlichen}}, Dissertation, University of T\"ubungen, Leipzig (1914).
\bibitem{Per}
O. Perron, {\it{\"Uber einen Satz des Herrn Poincar$\acute{\mbox{e}}$}}, J. Reine Angew. Math. {\bf{136}}, 17 (1909)
\bibitem{BA1}
C. R. Adams, {\it{On the irregular cases of linear ordinary difference equations}}, Trans. Amer. Math. Soc. {\bf{30}}, 507 (1928).
\bibitem{BA2}
G. D. Birkhoff, {\it{General theory of linear difference equations}}, Trans. Amer. Math. Soc. {\bf{12}}, 243 (1911).
\bibitem{BA3}
G. D. Birkhoff, {\it{Formal theory of irregular linear difference equations}}, Acta Math. {\bf{54}}, 205 (1930).
\bibitem{BA4}
G. D. Birkhoff and W. J. Trjitzinsky, {\it{Analytic theory of singular difference equations}}, Acta Math. {\bf{60}}, 1 (1932).
\bibitem{LI}
H. Li, {\it{Asymptotic expansions for second-order linear difference equations}}, J. Comput. Appl. Math. {\bf{41}}, 65 (1992).
\bibitem{Wong}
R. Wong and H. Li, {\it{Asymptotic expansions for second-order linear difference equations II}}, Stud. Appl. Math. {\bf{87}}, 289 (1992).
\end{thebibliography}
\end{document}